\documentclass{amsart}
\numberwithin{equation}{section}

\newtheorem{thm}{Theorem}[section]
\newtheorem{lem}[thm]{Lemma}

\newtheorem{defin}[thm]{Definition}

\newtheorem{remark}[thm]{Remark}

\begin{document}

\begin{center}
\textbf{{\large {\ A linear inverse problem with semi-nonlocal boundary conditions for a three-
dimensional equation of mixed type of the second kind of the fourth order in a parallelepiped. }}}\\[0pt]
\medskip \textbf{Sirojiddin Dzhamalov $^{1,2}$, Bakhtiyor Khalkhadzhaev $^{1}$}\\[0pt]
\textit{siroj63@mail.ru, xalxadjaev@yandex.ru\\[0pt]}
\medskip \textit{\ Institute of Mathematics, Academy of Science of Uzbekistan}$^{1}$

\medskip \textit{\ University of Tashkent for Applied Sciences}$^{2}$

\end{center}

\textbf{Abstract}: In~this article the~correctness of~a~linear inverse problem with~semi-nonlocal boundary conditions for~a~three-dimensional equation in~a~parallelepiped is~considered. The~equation itself is~a~fourth order  mixed type equation of~the~second kind. The~existence and uniqueness theorems for~a~generalized solution of~the~inverse problem in~a~certain class of~integrable functions are~proved using the~methods of~Fourier, "$\varepsilon $-regularization", a~priori estimates, approximating sequences and~contracting mappings.



\section{Introduction}
\quad In~the~process of~studying nonlocal problems, a~close relationship between problems with~nonlocal boundary conditions and inverse problems was revealed (see for example, \cite{L3},\cite{L8}-\cite{L11}). Inverse problems for~classical equations (parabolic, elliptic, and hyperbolic) are studied quite well (see for example, \cite{L1},\cite{L4},\cite{L17},\cite{L24},\cite{L25}). Linear inverse problems for~the~second order mixed type equations in~the~plane were studied by~A.G.~Megrabov, K.B.~Sabitov and their students. For~ multidimensional the~second order mixed type equations of~the~first and second kinds in~bounded domains were studied by S.Z. Dzhamalov, R.R. Ashurov, S.G. Pyatkov and A.I.Kozhanov\cite{L12}-\cite{L14}, \cite{L21},\cite{L23},\cite{L28}.

\quad Forward problems for~high-order mixed-type equations were studied in~\cite{L6},\cite{L7}, \cite{L15},\cite{L27},
\cite{L29} but inverse problems for~high-order mixed-type equations were practically not studied. The~aim of~this work is to fill this gap.

\quad In~this article, a~new method is proposed to~study the~unique solvability of~inverse problems for~a~ three-dimensional mixed-type fourth-order equation of~the~second kind in~a~rectangle. This method is based on~reducing the~inverse problem to forward semi-nonlocal boundary value problems for~a~family of~loaded fourth-order mixed type differential equations of~the~second kind.

\quad Recall that a~loaded equation is a~partial differential equation that contains the~values of~certain functionals from the~solution either in~the~coefficients or on~the~right-hand side \cite{L20},\cite{L25}.

\quad In~a~three-dimensional parallelepiped
$$G=\{ (x,y, t)\mid 0<x<1\, \, \, 0<y<\ell ;0<t<T<+\infty \} $$
we consider the fourth-order mixed-type equations of the second kind:
\begin{equation} \label{GrindEQ__1_}
Lu=Pu-Mu+Nu=\psi (x,t,y)\, .
\end{equation}

Here $Pu=\sum\limits_{i=0}^{\ 4}K_{i} (x,t) D_{t}^{i} u;$ $ Mu=au_{xxxx} +bu_{xxtt} -cu_{xx}; Nu=u_{yyyy} \, .$

\quad Let the~following conditions be satisfied for~the~coefficients of~the~equation \eqref{GrindEQ__1_}:

\quad $K_{4} (t)\in C^{3} (0,T)\cap C[0,T];\; K_{4} (0)=K_{4} (T)=0;\, K_{i} (x,t)\in C^{2} (Q)\cap C(\overline{Q})$ for all $x\in \left[0,1\right],$ $a,b,c$ are positive constants, $D_{t}^{i} u=\frac{\partial ^{i} u }{\partial t^{i} } \, \, ,\, $ $(i=1,2,3,4),$ $D_{t}^{0} u=u.$

\quad Equation \eqref{GrindEQ__1_} is~an~equation of~mixed type of~the~second kind, since there are no restrictions are imposed on~the~sign of~the~function $K_{4} (t)$ with respect to a~variable $t$ inside the~segment  $[0,T]$ \cite{L3}-\cite{L5},\cite{L9}.

\quad Hereinafter, we will assume that $\psi \, (x,t,y)=g\, (x,t,y)+h\, (x,t)\, f\, (x,t,y),$ where $g(x,t,y)$ and  $f(x,t,y)$ are given functions, and the~function $h\, (x,t)$ is to be determined.
\quad \textbf{\textit{Linear inverse problem:}} Find functions $\{ u(x,t,y),\; h(x,t)\} ,$ that satisfy equation \eqref{GrindEQ__1_} in~the~domain $G$ such that the~function $u(x,t,y)$ satisfies the following boundary conditions
\begin{equation} \label{GrindEQ__2_}
\gamma D_{t}^{p} \left. u\right|_{t=0} =D_{t}^{p} \left. u\right|_{t=T} ;\, \, \, p=0,1,2\,
\end{equation}
\begin{equation} \label{GrindEQ__3_}
\left. u\right|_{x=0} =\left. u\right|_{x=1} =0;
\end{equation}
\begin{equation} \label{GrindEQ__4_}
\left. u_{xx} \right|_{x=0} =\left. u_{xx} \right|_{x=1} =0\, \, \, \,
\end{equation}
\begin{equation} \label{GrindEQ__5_}
\left. u\right|_{y=0} =\left. u\right|_{y=1} =0;
\end{equation}
\begin{equation} \label{GrindEQ__6_}
\left. u_{yy} \right|_{y=0} =\left. u_{yy} \right|_{y=\ell } =0\, \, \, \,
\end{equation}
and auxiliary condition
\begin{equation} \label{GrindEQ__7_}
u(x,t,\ell _{0} )=\varphi _{0} \, (x,t), \text{ where }  0<\ell _{0} <\ell <+\infty
\end{equation}
and the~function $h(x,t)$ belongs to~the~following class
$$U=\{ \left. (u,h)\right|\, u\in W_{2}^{4,3} (G);h\in W_{2}^{4} (Q)\} .$$
Here  $W_{2}^{l,s} (G)$ denotes the space of functions with the norm
$$ \left\| u\right\| _{W_{2}^{l,s} (G)}^{2} =\sqrt{\frac{2}{\ell } } \sum _{k=1}^{\infty }(1+\mu _{k}^{4} )^{s} \left\| u_{k} \, (x,t)\right\| _{W_{2}^{l} (Q)}^{2},     \eqno(A) $$
where $\; s\ge 3, l\geq4 $ are positive integers, and the norm in the Sobolev space $W_{2}^{l} (Q)$ is defined as:
$$\left\| \vartheta \right\| _{W_{2}^{l} (Q)}^{2} =\sum\limits_{\left|\alpha \right|\le l}\, \int _{Q}\left|D^{\alpha } \vartheta \right|^{2}   dxdt,$$
where $\alpha $ is a multi-index, $D^{\alpha } $ is a generalized derivative with respect to variables $x$ and $t$.
Finally, $u_{k}(x,t)$ denotes the Fourier coefficient of the function $u(x,t,y),$ i.e.,
$$u_{k}(x,t)=\sqrt{\frac{2}{\ell } } \int_{0}^{\ell }u(x,t,y)\sin\mu_{k} ydy$$
Obviously,   $W_{2}^{l,s} (G)$ with norm (A) is a Banach space \cite{L11}-\cite{L14}.

\begin{defin}\label{def1}A regular solution to the inverse problem \eqref{GrindEQ__1_}-\eqref{GrindEQ__7_} is a function $\{u,h\}\in U, $ that satisfies equation \eqref{GrindEQ__1_} almost everywhere in the domain $G$ with conditions \eqref{GrindEQ__2_}-\eqref{GrindEQ__7_}.
\end{defin}

We will prove the unique solvability of the problem \eqref{GrindEQ__1_}-\eqref{GrindEQ__7_} using the Fourier method, i.e. to find a solution to problem \eqref{GrindEQ__1_}-\eqref{GrindEQ__7_}, we apply the Fourier method in variable  $y$ to the problem.

\quad The solution to problem \eqref{GrindEQ__1_}-\eqref{GrindEQ__7_} is sought in the following form:
$$ u(x,t,y)=\sum_{k=1}^{\infty}u_{k} (x,t)Y_{k}(y)  \eqno(B) $$
where the functions $\, Y_{k} (y)=\{ \sqrt{\frac{2}{\ell } } \sin \mu _{k} y\} $,$\mu _{k} =\frac{\pi k}{\ell } ,$ $k=1,2,3,...$ are solutions to the spectral problem of a fourth-order equation with Dirichlet conditions, that is
$$\begin{array}{l} {Y_{k}^{(IV)} =\mu _{k}^{4} Y_{k} } \\ {Y_{k} (0)=Y_{k} (\ell )=0} \\ {Y_{k} ^{{'} {'} } (0)=Y_{k} ^{{'} {'} } (\ell )=0\, .} \end{array}$$
It is known that the system of eigenfunctions $\{ Y_{k} (y)\}$ is fundamental in space $L_{2} [\, 0,\ell \; ]$ and forms an orthonormal basis in it \cite{L16},\cite{L22}, and functions $u_{k} (x,t)\, \; k=\, 1,\, 2,\, 3,...$ are unknown coefficients.We assume the following conditions to be fulfilled.

Let all the coefficients of equation \eqref{GrindEQ__1_} be sufficiently smooth functions in the domain $\overline{Q}$ and let the following conditions be satisfied regarding the coefficients of the right-hand side and the given function $\varphi _{0} (x,t)$.

\textbf{Condition 1:}

periodicity $K_{4t} (0)=K_{4t} (T);\; K_{i} (x,0)=K_{i} (x,T);i=0,2,3,$ for every $ x\in \left[0,1\right].$

non-locality:          $\gamma g(x,0,z)=g(x,T,z),$$\gamma f(x,0,z)=f(x,T,z).$

smoothness:  $f(x,t,l_{0} )=f_{0} (x,t)\in {C_{x,t}^{0,1} (Q)},$ $\left|f_{0}(x,t)\right|\ge\eta >0;$ $f\in W_{2}^{3,3} (G),$ $ g\in W_{2}^{1,3} (G).$

\textbf{Condition 2:}
$$\varphi _{0} (x,t)\in W_{2}^{5} (Q); \, \gamma D_{t}^{p} \left. \varphi _{0} \right|_{t=0} =\left. D_{t}^{p} \varphi _{0} \right|_{t=T} ,p=0,1,2,3; $$
$$\left. \varphi _{0} \right|_{x=0} =\left. \varphi _{0} \right|_{x=1} =0,\; \left.\varphi _{0xx}\right|_{x=0} =\left.\varphi _{0xx} \right|_{x=1} =0.$$

To formulate the main result, it is necessary to perform some formal constructions.

We consider the traces of equation \eqref{GrindEQ__1_} at $y=\ell _{0}$, i.e.
$$\left. Lu\right|_{y=\ell _{0} } =P\varphi _{0} -M\varphi _{0} +\sqrt{\frac{2}{\ell } } \sum _{k=1}^{\infty }\mu _{k}^{4} u_{k}  \sin \mu _{k} \ell _{0} =g_{0} (x,t)\, \, +h(x,t)f_{0} (x,t)$$
Now, considering the auxiliary condition \eqref{GrindEQ__7_} and that $f_0\neq 0$, we define the unknown function $h(x,t)$ in the form of an integral
$$ h(x,t)=\frac{1}{f_{0} (x,t)} \left[\Phi _{0} +\sqrt{\frac{2}{\ell } } \sum _{k=1}^{\infty }\mu _{k}^{4}  u_{k} (x,t)\sin \mu _{k} \ell _{0} \right]$$
where$\; \Phi _{0} =L_{0} \varphi _{0} -g_{0};$ $ g(x,t,\ell _{0} )=g_{0} (x,t);$  $L_{0} \varphi _{0} =P\varphi _{0} -M\varphi _{0};$
$P\varphi _{0} =\sum\limits_{i=0}^{4}k_{i} (x,t)D_{t}^{i} \varphi _{i0};$ $M\varphi _{0} =aD_{x}^{4} \varphi _{0} +bD_{x,t}^{2,2} \varphi _{0} -cD_{x}^{2} \varphi _{0}$ and to determine the functions in domain $u_{k} (x,t)$ $Q=(0,1)\times (0,T)$ we obtain countably many loaded systems of equations of mixed type of the fourth order:
\begin{equation} \label{GrindEQ__8_}
\begin{array}{c} {\Im u_{k} =Pu_{k} -Mu_{k}+\mu_{k}^{4}u_{k}=g_{k} (x,t)+} \\ {+\frac{f_{k}(x,t)}{f_{0}(x,t)}[\Phi_{0}+\sqrt{\frac{2}{\ell}} \sum\limits_{m=1}^{\infty}\mu_{m}^{4}u_{m} (x,t)\sin\mu_{m}\ell_{0}]=F(u_{k})\quad}\end{array}
\end{equation}
with semi-nonlocal boundary conditions:
\begin{equation} \label{GrindEQ__9_}
\gamma D_{t}^{p} \left. u_{k} \right|_{t=0} =\left. D_{t}^{p} u_{k} \right|_{t=T} ;p=0,1,2,
\end{equation}
\begin{equation} \label{GrindEQ__10_}
\left. u_{k} \right|_{x=0} =\left. u_{k} \right|_{x=1} =0,
\end{equation}
\begin{equation} \label{GrindEQ__11_}
\left. u_{kxx} \right|_{x=0} =\left. u_{kxx} \right|_{x=1} =0\, \,
\end{equation}
The main results are the following two theorems:

\begin{thm}\label{Siroj1}
Let the conditions 1 and 2 be satisfied and let \\ $ (2K_{3} +(2j-3)K_{4t} +3\lambda K_{4} )\ge \delta _{3} >0,\; j=0,1,2;$ $-2K_{1} +K_{2t} -\lambda K_{2} \ge \delta _{2} >0,$ $-\lambda K_{0} +K_{0t} \ge \delta _{1} >0,$ for any $ (x,t)\in \mathop{\overline{Q}}\limits^{} ,$ where $ \lambda =\frac{2}{T} \ln \left|\gamma \right|>0,\, \, \, \, \left|\gamma \right|>1\, ,$ and let there exist a positive number $\sigma $ such that estimates $ \delta _{0} -31e^{\lambda T} \sigma ^{-1} >\delta >0,$ $ q=M\left\| f\right\| _{W_{2}^{3,3} (G)}^{2} <1$ are hold for $\, \delta _{0} =\min \{\delta _{1} ,\lambda\{a,b,c\} ,\delta_{2},\delta_{3}\} $,  where $ M=\sigma \, e^{\lambda T} \lambda ^{4} m\, \delta ^{-1} \ell^{-2} \eta ^{-2} \left\| f_{0} \right\| ^{2} _{C_{x,t}^{0,1} (Q)}$,  $m=10c_{1} c_{2} c_{3}$, \\$c_{1} =\sum\limits_{k=1}^{\infty}\frac{\mu_{k}^{8}}{(1+\mu_{k}^{4})^{3}}, \,\, c_{i}(i=2,3)$ are the coefficients of the Sobolev embedding theorem. Then there is a unique solution to the linear inverse problem \eqref{GrindEQ__1_}-\eqref{GrindEQ__7_} from the indicated class $U.$
\end{thm}


Proof. \textit{Theorem \eqref{Siroj1} is proven according to the following scheme. }

1. Let us show that the solution $\{u,h\}\in U$ to the problem \eqref{GrindEQ__1_}-\eqref{GrindEQ__6_} satisfies auxiliary condition \eqref{GrindEQ__7_}.\textit{}

2. To prove the unique solvability of problem \eqref{GrindEQ__8_}-\eqref{GrindEQ__11_}, first we investegate the unique solvability of the auxiliary problem \eqref{GrindEQ__17_}-\eqref{GrindEQ__20_}, that is we investegate the solvability of a family of nonlinear loaded composite type differential equations of with a small parameter.

3. Then, using this auxiliary problem \eqref{GrindEQ__17_}-\eqref{GrindEQ__20_}, we investegate the unique solvability of a family of nonlinear loaded mixed-type fourth-order equations of the second kind \eqref{GrindEQ__8_}-\eqref{GrindEQ__11_}.

4. Using the unique solvability of problem \eqref{GrindEQ__8_}-\eqref{GrindEQ__11_}, we prove the unique solvability of the linear inverse problem \eqref{GrindEQ__1_}-\eqref{GrindEQ__7_}.

Now let us start realise this scheme.
To obtaining various a priori estimates, we use Cauchy's inequality with $\sigma $  that is
$$\forall u,\vartheta \ge 0;\; \forall \sigma >0;\; \quad 2u\cdot \vartheta \le \sigma u^{2} +\sigma ^{-1} \vartheta ^{2} .$$

Let us show that the function $u(x,t,y)$ defined by series (B) satisfies the auxiliary condition \eqref{GrindEQ__7_}, i.e. $ u(x,t,\ell _{0} )=\varphi _{0} (x,t).$

Assuming the opposite $u(x,t,\ell _{0} )\ne \varphi _{0} (x,t),$ for  function$ $ $z(x,t)=u(x,t,\ell _{0} )-\varphi _{0} (x,t)=\sum \limits_{k=1}^{\infty }u_{k} (x,t) \sin \mu _{k} \ell _{0} -\varphi _{0} (x,t)$ in the domain $ Q $ from \eqref{GrindEQ__8_}-\eqref{GrindEQ__11_}, multiplying \eqref{GrindEQ__8_}-\eqref{GrindEQ__11_} by $\sin \mu _{s} \ell _{0} $ and summing over $k$ from $ 1 $ to $\infty ,$ the following identity is obtained:
\begin{equation} \label{GrindEQ__12_}
\begin{array}{l} {\sum\limits_{k=1}^{\infty }\Im u_{k}\sin\mu_{k}\ell_{0}=P\sum\limits_{k=1}^{\infty }u_{k}\sin\mu_{k}\ell_{0}-M\sum\limits _{k=1}^{\infty }u_{k}\sin\mu_{k}\ell_{0}
 +\sqrt{\frac{2}{\ell } }\sum\limits_{k=1}^{\infty }\mu_{k}^{4}u_{k}\sin\mu_{k}\ell_{0}}= \\
 {=\sum\limits_{k=1}^{\infty }g_{k}\sin\mu_{k}\ell_{0}+\frac{\sum\limits_{k=1}^{\infty}f_{k}\sin\mu_{k}\ell_{0}}{f_{0} (x,t)}\, [\Phi_{0} +\sqrt{\frac{2}{\ell} }\sum\limits_{m=1}^{\infty }\,\mu_{m}^{4}\,u_{m}\sin\mu_{m}\ell_{0} ].} \end{array}
\end{equation}

Considering boundary conditions \eqref{GrindEQ__9_}--\eqref{GrindEQ__11_}, \eqref{GrindEQ__12_} we obtain the following problem:
\begin{equation} \label{GrindEQ__13_}
L_{0} \, z=\sum _{i=0}^{4}K_{i} \,  D\, _{t}^{i} \, z-M\, z=0
\end{equation}
with semi-nonlocal boundary conditions:
\begin{equation} \label{GrindEQ__14_}
\gamma D_{t}^{p} \left. z\right|_{t=0} =D_{t}^{p} \left. z\right|_{t=T} ;\, \, \, p=0,1,2,
\end{equation}
\begin{equation} \label{GrindEQ__15_}
\left. z\right|_{x=0} =\left. z\right|_{x=1} =0;
\end{equation}
\begin{equation} \label{GrindEQ__16_}
\left. z_{xx} \right|_{x=0} =\left. z_{xx} \right|_{x=1} =0\, \, \, .
\end{equation}
Now we prove the uniqueness of the solution to problem \eqref{GrindEQ__13_}, \eqref{GrindEQ__16_}.

To do this, consider the identity $2(L_{0} z,e^{-\lambda t} z_{t} )_{0} =0.$ Integrating it by parts, considering conditions of Theorem \eqref{Siroj1} and conditions \eqref{GrindEQ__14_}--\eqref{GrindEQ__16_}, we obtain an inequality $\left\| z\right\| _{2} \le 0,$ which implies that $z(x,t)=0,$ and that problem \eqref{GrindEQ__13_}-\eqref{GrindEQ__16_} has a unique solution. This means that the solution to problem \eqref{GrindEQ__1_}-\eqref{GrindEQ__6_} satisfies the additional condition \eqref{GrindEQ__7_}, i.e. $u(x,t,\ell _{0} )=\varphi _{0} (x,t)$.
\section{A family of loaded differential equations of the fifth order with a small parameter(auxiliary problem)}

We will prove the solvability of problem \eqref{GrindEQ__8_}-\eqref{GrindEQ__11_} using the methods of ``$\varepsilon-$regularization'', successive approximations and a priori estimates, namely, we will consider in domain $Q=(0,1)\times (0,T)$ a semi-nonlocal boundary value problem for a family of loaded differential equations of the fifth order with a small parameter:
\begin{equation} \label{GrindEQ__17_}
\begin{array}{l} {\Im _{\varepsilon } u_{k} =-\varepsilon \frac{\partial \, \Delta ^{2} u_{k,\varepsilon } }{\partial \, t} +Pu_{k,\varepsilon } -Mu_{k,\varepsilon } +\mu _{k}^{4} u_{k,\varepsilon } =g_{k} (x,t)+} \\ {+\frac{f_{k} \, (x,t)}{f_{0} (x,t)} [\Phi _{0} +\sqrt{\frac{2}{\ell } } \sum\limits_{m=1}^{\infty }\mu _{m}^{4}  u_{m,\varepsilon } (x,t)\sin\mu_{m}\ell_{0}  ]=F(u_{k,\varepsilon } ),\quad } \end{array}
\end{equation}
with semi-nonlocal boundary conditions:
\begin{equation} \label{GrindEQ__18_}
\gamma D_{t}^{p} \left. u_{k,\varepsilon } \right|_{t=0} =\left. D_{t}^{p} u_{k,\varepsilon } \right|_{t=T} ;p=0,1,2,3,4,
\end{equation}
\begin{equation} \label{GrindEQ__19_}
\left. u_{k,\varepsilon } \right|_{x=0} =\left. u_{k,\varepsilon } \right|_{x=1} =0,
\end{equation}
\begin{equation} \label{GrindEQ__20_}
\left. u_{k,\varepsilon xx} \right|_{x=0} =\left. u_{k,\varepsilon xx} \right|_{x=1} =0\, ,
\end{equation}
where $\varepsilon $ is a small positive number.

When proving Theorem~\ref{Siroj1} and the well-posedness of problem \eqref{GrindEQ__17_}-\eqref{GrindEQ__20_}, we will use the following notation and auxiliary lemmas.
\\
Let us define spaces of vector functions $W_{i} (Q)=\{ \vartheta _{k} :\vartheta _{k} \in W_{2}^{i} (Q)\} $     where $i=0,1,2,3,4;\; \; k=1,2,3,...$  with finite norm

$$ \left\langle \vartheta _{k} \right\rangle _{i}^{2} =\sum _{k=1}^{\infty }(1+\mu _{k}^{4} )^{3}  \left\| \vartheta _{k} \right\| _{W_{2}^{i} (Q)}^{2} ; i=0,1,2,3,4. \eqno (C) $$

Space  $W_{i} (Q)\; $ with the above norm (C) is a Banach space.

From the definition of spaces $W\; _{2}^{i} (Q),\, \, i=0,1,2,3,4$  \quad it follows
$$W_{4} (Q)\subset W_{3} (Q)\subset W_{2} (Q)\subset W_{1} (Q)\subset W_{0} (Q).$$

Below $W(Q)$ will denote the class of vector functions $\{ \vartheta _{k} (x,t)\} _{k=1}^{\infty } $ such that  $\{ \vartheta _{k} (x,t)\} _{k=1}^{\infty } \in W_{4} (Q)$, $\{ \frac{\partial }{\partial t} \Delta ^{2} \vartheta _{k} \} _{1}^{\infty } \, \in W_{0} (Q)$ and satisfying the boundary conditions \eqref{GrindEQ__18_}, \eqref{GrindEQ__20_}.

\begin{defin}\label{def2} The solution to problem \eqref{GrindEQ__17_}-\eqref{GrindEQ__20_} is called a vector function $\{ u_{s,\varepsilon } (x,t)\} \in W$ that satisfies equation \eqref{GrindEQ__17_} almost everywhere in the domain $Q$.
\end{defin}

Now we will prove the solvability of problem \eqref{GrindEQ__17_}-\eqref{GrindEQ__20_} in the domain $Q$ using the methods of successive approximations \cite{L2},\cite{L22}.
\begin{equation} \label{GrindEQ__21_}
\begin{array}{l} {\Im _{\varepsilon } u_{_{k,\varepsilon } }^{(l)} =-\varepsilon \frac{\partial \, \Delta ^{2} u_{k,\varepsilon }^{(l)} }{\partial \, t} +Pu_{k,\varepsilon }^{(l)} -Mu_{k,\varepsilon }^{(l)} +\mu _{k}^{4} u_{k,\varepsilon }^{(l)} =g_{k} (x,t)+} \\ {+\frac{f_{k} \, (x,t)}{f_{0} (x,t)} [\Phi _{0} +\sqrt{\frac{2}{\ell } } \sum\limits_{m=1}^{\infty }\mu _{m}^{4}  u_{m,\varepsilon }^{(l)} (x,t)\sin \mu _{m} \ell _{0}  ]=F(u_{k,\varepsilon }^{(l-1)} )\quad } \end{array}
\end{equation}
with semi-nonlocal boundary conditions:
\begin{equation} \label{GrindEQ__22_}
\gamma D_{t}^{p} \left. u_{k,\varepsilon }^{(l)} \right|_{t=0} =\left. D_{t}^{p} u_{k,\varepsilon }^{(l)} \right|_{t=T} ;p=0,1,2,3,4,
\end{equation}
\begin{equation} \label{GrindEQ__23_}
\left. u_{k,\varepsilon }^{(l)} \right|_{x=0} =\left. u_{k,\varepsilon }^{(l)} \right|_{x=1} =0,
\end{equation}
\begin{equation} \label{GrindEQ__24_}
\left. u_{k,\varepsilon xx}^{(l)} \right|_{x=0} =\left. u_{k,\varepsilon xx}^{(l)} \right|_{x=1} =0\, ,
\end{equation}
where $\varepsilon >0,\quad l=0,1,2,...\quad ,\left\{u_{\varepsilon }^{\left(-1\right)} \right\}=0\, .$

\begin{lem} \label{Sir1}
Let all the conditions of Theorem \eqref{Siroj1} be satisfied. Then for the solution of  problem \eqref{GrindEQ__21_}-\eqref{GrindEQ__24_} the following estimates are valid
$$I).\quad \frac{\varepsilon }{\delta } (\left\langle u_{k,\varepsilon ttt}^{(l)} \right\rangle _{0}^{2} +\left\langle u_{k,\varepsilon ttx}^{(l)} \right\rangle _{0}^{2} +\left\langle u_{k,\varepsilon ttt}^{(l)} \right\rangle _{0}^{2} )+\left\langle u_{k,\varepsilon }^{(l)} \right\rangle _{2}^{2} \le const(\not{l}, \not{\varepsilon} ,\not{\mu_{k}}\ ),$$
$$II).\quad \; \frac{\varepsilon }{\delta } \left\langle \frac{\partial \Delta ^{2} u_{k,\varepsilon }^{(l)} }{\partial t} \right\rangle _{0}^{2} +\left\langle u_{k,\varepsilon }^{(l)} \right\rangle _{4}^{2} \le const(\not{l}, \not{\varepsilon} ,\not{\mu_{k}}\ ).\quad\quad\quad\quad\quad\quad\quad\quad\quad\quad$$
Here and below symbol $const(\not{l}, \not{\varepsilon} ,\not{\mu_{k}}\ )$ denotes a constant that does not depend on parameters  $l,\varepsilon ,\mu _{k}  $.
\end{lem}

Proof of Lemma \ref{Sir1}. Consider the following identity:
\begin{equation} \label{GrindEQ__25_}
\left|-2\left(\Im _{\varepsilon } u_{k,\varepsilon }^{(l)} ,e^{-\lambda t} u_{k,\varepsilon t}^{(l)} \right)_{0}
\right|=\left|-2\left(F\left(u_{k,\varepsilon }^{(l-1)} \right),e^{-\lambda t} u_{k,\varepsilon t}^{(l)} \right)_{0} \right|
\end{equation}
where constants $\lambda >0$ will be choosen  later.

Due to the conditions of Theorem \eqref{Siroj1}, boundary conditions \eqref{GrindEQ__22_}-\eqref{GrindEQ__24_}, considering that $\gamma ^{2} =e^{\lambda T}$, by integrating identity \eqref{GrindEQ__25_} by parts and applying Cauchy's inequalities with $\sigma $ of [12], from identity \eqref{GrindEQ__25_} we obtain the first recurrent formula
\begin{equation} \label{GrindEQ__32_}
\begin{array}{c} {\frac{\varepsilon }{\delta } (\left\langle u_{k,\varepsilon ttt}^{(l)} \right\rangle _{0}^{2} +\left\langle u_{k,\varepsilon ttx}^{(l)} \right\rangle _{0}^{2} +\left\langle u_{k,\varepsilon txx}^{(l)} \right\rangle _{0}^{2} )+\left\langle u_{k,\varepsilon }^{(l)} \right\rangle _{2}^{2} \le } \\ {+2\sigma \delta ^{-1} e^{\lambda T} \left[\left\langle g_{k} \right\rangle _{0}^{2} +2\eta ^{-2} c_{2} \left\langle f_{k} \right\rangle _{2}^{2} \left(T_{0}^{2} \left\| \varphi _{0} \right\| _{W_{2}^{4} (Q)}^{2} +\left\| g_{0} \right\| _{0}^{2} \right)\right]+} \\ {+2c_{1} c_{2} \sigma \delta ^{-1} \ell ^{-2} \eta ^{-2} e^{\lambda T} \left\langle f_{k} \right\rangle _{2}^{2} \left\langle u_{k,\varepsilon }^{(l-1)} \right\rangle _{2}^{2} .} \end{array}
\end{equation}

We introduce the notation
$$ 2\sigma \delta ^{-1} e^{\lambda T} \left[\left\langle g_{k} \right\rangle _{0}^{2} +2\eta ^{-2} c_{2} \left\langle f_{k} \right\rangle _{2}^{2} \left(T_{0}^{2} \left\| \varphi _{0} \right\| _{W_{2}^{4} (Q)}^{2} +\left\| g_{0} \right\| _{0}^{2} \right)\right]=A_{1} $$ and using  conditions of Theorem \eqref{Siroj1}
$$2c_{1} c_{2} \sigma \delta ^{-1} \ell ^{-2} \eta ^{-2} e^{\lambda T} \left\langle f_{k} \right\rangle _{2}^{2} <q=M\left\langle f_{k} \right\rangle _{3}^{2} <1\, ,$$ from the recurrent formula \eqref{GrindEQ__32_}, we obtain the validity of the estimate I). For this, we take function $\left\{{u}_{\varepsilon }^{(-1)} \right\}\equiv \left\{0\right\} $ as an ``initial approximation''. Hence
$$\begin{array}{c} {\frac{\varepsilon }{\delta } (\left\langle u_{k,\varepsilon ttt}^{(0)} \right\rangle _{0}^{2} +\left\langle u_{k,\varepsilon ttx}^{(0)} \right\rangle _{0}^{2} +\left\langle u_{k,\varepsilon txx}^{(0)} \right\rangle _{0}^{2} )+\left\langle u_{k,\varepsilon }^{(0)} \right\rangle _{2}^{2} \le } \\ {\le 2\sigma \delta ^{-1} e^{\lambda T} \left[\left\langle g_{k} \right\rangle _{0}^{2} +2\eta ^{-2} c_{2} \left\langle f_{k} \right\rangle _{2}^{2} \left(T_{0}^{2} \left\| \varphi _{0} \right\| _{W_{2}^{4} (Q)}^{2} +\left\| g_{0} \right\| _{0}^{2} \right)\right]\equiv A_{1}.} \end{array}$$

Continuing this process, by induction we obtain the first a priori estimate for any function $\left\{u_{k,\varepsilon }^{(l)} \right\},\; \forall l\ge 0$
\begin{equation} \label{GrindEQ__33_}
\frac{\varepsilon }{\delta } (\left\langle u_{k,\varepsilon ttt}^{(l)} \right\rangle _{0}^{2} +\left\langle u_{k,\varepsilon ttx}^{(l)} \right\rangle _{0}^{2} +\left\langle u_{k,\varepsilon txx}^{(l)} \right\rangle _{0}^{2} )+\left\langle u_{k,\varepsilon }^{(l)} \right\rangle _{2}^{2} \le A_{1}\cdot \sum _{k=0}^{l}q^{k}  <\frac{A_{1}}{1-q}
\end{equation}
Now we prove the validity of estimate II ).
To do this, consider the following identity:
\begin{equation} \label{GrindEQ__34_}
\left|-2\int\limits _{Q}\Im _{\varepsilon } u_{k,\varepsilon }^{(l)} \; e^{-\lambda t} Pu_{k,\varepsilon }^{(l)} \, dxdt  \right|= \left|-2\int \limits _{Q}F\left(u_{k,\varepsilon }^{(l-1)} \right) e^{-\lambda t} Pu_{k,\varepsilon }^{(l)} dxdt  \right|\, ,\quad
\end{equation}
where $Pu_{k,\varepsilon }^{(l)} \equiv \frac{\partial \Delta ^{2} u_{k,\varepsilon }^{(l)} }{\partial \; t} -2\lambda \frac{\partial ^{2} \Delta u_{k,\varepsilon }^{(l)} }{\partial \; t^{2} } +3\lambda ^{2} \frac{\partial \Delta u_{k,\varepsilon }^{(l)} }{\partial t} -\frac{\lambda u_{k,\varepsilon tt}^{(l)} }{2} +\frac{\lambda ^{2} u_{k,\varepsilon }^{(l)} }{16}. $
$$\frac{\partial\Delta^{2}u_{k,\varepsilon}^{(l)}}{\partial t}=\frac{\partial}{\partial t}\left(u_{k,\varepsilon tttt}^{(l)}+2u_{k,\varepsilon ttxx}^{(l)}+u_{k,\varepsilon xxxx}^{(l)}\right)$$

Integrating \eqref{GrindEQ__34_} by parts, taking into account the conditions of Theorem \eqref{Siroj1} and boundary conditions \eqref{GrindEQ__22_}, \eqref{GrindEQ__24_} we obtain the second recurrent formula
\begin{equation} \label{GrindEQ__40_}
\begin{array}{l} {\frac{\varepsilon }{\delta } \left\langle \frac{\partial \Delta ^{2} u_{k,\varepsilon }^{(l)} }{\partial t} \right\rangle _{0}^{2} +\left\langle u_{k,\varepsilon }^{(l)} \right\rangle _{4}^{2} \quad \le } \\ {\le 15\sigma \lambda ^{4} e^{\lambda T} \cdot [\left\langle g_{k} \right\rangle _{1}^{2} +2\eta ^{-2} c_{2} \left\| f_{0} \right\| ^{2} _{C_{x,t}^{0,1} (Q)} \left\langle f_{k} \right\rangle _{3}^{2} [T_{1}^{2} \left\| \varphi _{0} \right\| _{5}^{2} +\left\| g_{0} \right\| _{1}^{2} ]]+} \\ {+6\sigma \lambda ^{4} \ell ^{-2} \eta ^{-2} e^{\lambda T} c_{1} c_{2} \left\| f_{0} \right\| ^{2} _{C_{x,t}^{0,1} (Q)} \left\langle f_{k} \right\rangle _{3}^{2} \left\langle u_{k,\varepsilon }^{(l-1)} \right\rangle _{4}^{2} .\quad \quad \quad } \end{array}
\end{equation}

Let us introduce the notation
\[15\sigma \lambda ^{4} e^{\lambda T} \cdot [\left\langle g_{k} \right\rangle _{1}^{2} +2\eta ^{-2} c_{2} \left\| f_{0} \right\| ^{2} _{C_{x,t}^{0,1} (Q)} \left\langle f_{k} \right\rangle _{3}^{2} [T_{1}^{2} \left\| \varphi _{0} \right\| _{5}^{2} +\left\| g_{0} \right\| _{1}^{2} ]]\equiv A_{2}
\]
and considering conditions of Theorem \eqref{Siroj1} i.e. $6\sigma \lambda ^{4} \ell ^{-2} \eta ^{-2} e^{\lambda T} c_{1} c_{2} \left\| f_{0} \right\| ^{2} _{C_{x,t}^{0,1} (Q)} \left\langle f_{k} \right\rangle _{3}^{2} <q=M\left\langle f_{k} \right\rangle _{3}^{2} <1$, from the recurrent formula \eqref{GrindEQ__40_}, we obtain the validity of estimate II). For this, we take the function as the ``initial approximation''. $\left\{u_{k,\varepsilon }^{(-1)} \right\}\equiv \left\{0\right\}.$
Then for the ``zero approximation'' we have
$$\begin{array}{l}\frac{\varepsilon}{\delta}\left\langle\frac{\partial\Delta^{2}u_{k,\varepsilon}^{(0)}}{\partial t}\right\rangle _{0}^{2}+\left\langle u_{k,\varepsilon }^{(0)}\right\rangle_{2}^{2}\le 15\sigma \lambda^{4}e^{\lambda T}[\left\langle g_{k}\right\rangle _{1}^{2}+\\
+2\eta^{-2}c_{2}\left\| f_{0} \right\| ^{2} _{C_{x,t}^{0,1} (Q)}\left\langle f_{k}\right\rangle_{3}^{2} [T_{1}^{2}\left\|\varphi_{0}\right\| _{5}^{2}+\left\|g_{0}\right\|_{1}^{2}]]\equiv A_{2}\end{array}$$
Continuing this process, by induction we obtain the second a priori estimate for any function$\left\{u_{k,\varepsilon }^{(l)} \right\},\; \forall l\ge 0$
$$\frac{\varepsilon }{\delta } \left\langle \frac{\partial \Delta ^{2} u_{k,\varepsilon }^{(l)} }{\partial t} \right\rangle _{0}^{2} +\left\langle u_{k,\varepsilon }^{(l)} \right\rangle _{2}^{2} \le A_{2} \cdot \sum _{k=0}^{l}q^{k}  \le \frac{A_{2} }{1-q} .$$
Similar the proof of estimate I ), estimate II ) is easily obtained.

The~Lemma \ref{Sir1} is proved.

Let's introduce a new function from $W(Q,{\mathbb R})$ the formula
$$\vartheta _{k,\varepsilon }^{(l)} =u_{k,\varepsilon }^{(l)} -u_{k,\varepsilon }^{(l-1)} \, ;\varepsilon >0:\, \, l=0,1,2,3...\, $$

Then the following lemma is valid.

\begin{lem} \label{Sir2}
Let all the conditions of Theorem \eqref{Siroj1} be satisfied. Then the following estimates hold for the function $\left\{\, \vartheta _{k,\varepsilon }^{(l)} \, \right\}\in W(Q,{\mathbb R}):$
$$III).\quad \frac{\varepsilon }{\delta } (\left\langle u_{k,\varepsilon ttt}^{(l)} \right\rangle _{0}^{2} +\left\langle u_{k,\varepsilon ttx}^{(l)} \right\rangle _{0}^{2} +\left\langle u_{k,\varepsilon txx}^{(l)} \right\rangle _{0}^{2} )+\left\langle u_{k,\varepsilon }^{(l)} \right\rangle _{2}^{2} \le A_{1}q^{(l)} ,$$
$$IV).\quad \; \frac{\varepsilon }{\delta } \left\langle \frac{\partial \Delta ^{2} u_{k,\varepsilon }^{(l)} }{\partial t} \right\rangle _{0}^{2} +\left\langle u_{k,\varepsilon }^{(l)} \right\rangle _{4}^{2} \le A_{2} q^{(l)}. \quad\quad\quad \quad\quad\quad\quad\quad\quad\quad$$
\end{lem}

Proof of Lemma \eqref{Sir2}. From \eqref{GrindEQ__21_}-\eqref{GrindEQ__24_} for the function $\left\{\, \vartheta _{k,\varepsilon }^{(l)} \, \right\}\in W(Q,{\mathbb R})$ \, we obtain the following problem
\begin{equation} \label{GrindEQ___41_}
\begin{array}{c}\Im _{\varepsilon } \vartheta _{k,\varepsilon }^{(l)} =-\varepsilon \frac{\partial \Delta ^{2} \vartheta _{k,\varepsilon }^{(l)} }{\partial t} +L_{0} \vartheta _{k,\varepsilon }^{(l)} +\mu _{k}^{4} \vartheta _{k,\varepsilon }^{(l)} =\\
=\sqrt{\frac{2}{\ell } } \frac{f_{k} \, (x,t)}{f_{0} (x,t)} \sum \limits _{k=1}^{\infty }\, \mu _{k}^{4} \vartheta _{k,\varepsilon }^{(l-1)} \sin \mu _{k}^{} \ell _{0}  \, \equiv F(\vartheta _{k,\varepsilon }^{(l-1)} )
\end{array}
\end{equation}
with semi-nonlocal boundary conditions\textbf{}
\begin{equation} \label{GrindEQ__42_}
\gamma D_{t}^{q} \left. \vartheta _{k,\varepsilon }^{(l)} \right|_{t=0} =D_{t}^{q} \left. \vartheta _{k,\varepsilon }^{(l)} \right|_{t=T} ;q=0,1,2,3,4,
\end{equation}
\begin{equation} \label{GrindEQ__43_}
\left. \vartheta _{k,\varepsilon }^{(l)} \right|_{x=0} =\left. \vartheta _{k,\varepsilon }^{(l)} \right|_{x=1} =0,
\end{equation}
where $\varepsilon >0,\quad l=0,1,2,...$;
Therefore, similar to the proof of Lemma \eqref{Sir1}, for function $\left\{\, \vartheta _{k,\varepsilon }^{(l)} \, \right\}=\left\{\, u_{k,\varepsilon }^{(l)} \, \right\}-\left\{\, u_{k,\varepsilon }^{(l-1)} \, \right\}\in W(Q,R)$we obtain the third recurrent formula
\begin{equation} \label{GrindEQ__44_}
\frac{\varepsilon }{\delta } (\left\langle \vartheta _{k,\varepsilon ttt}^{(l)} \right\rangle _{0}^{2} +\left\langle \vartheta _{k,\varepsilon ttx}^{(l)} \right\rangle _{0}^{2} +\left\langle \vartheta _{k,\varepsilon txx}^{(l)} \right\rangle _{0}^{2} )+\left\langle \vartheta _{k,\varepsilon }^{(l)} \right\rangle _{2}^{2} \le q\left\langle \vartheta _{k,\varepsilon }^{(l-1)} \right\rangle _{2}^{2}
\end{equation}
repeating the reasoning of Lemma \eqref{Sir1}, from \eqref{GrindEQ__44_} we obtain an a priori estimate III).
Estimate IV )  is proven in a similarly.

Lemma \ref{Sir2} is proven.

\begin{thm}\label{Siroj2}
Let all the conditions of Theorem \eqref{Siroj1} be satisfied. Then problem \eqref{GrindEQ__21_}-\eqref{GrindEQ__24_} is uniquely solvable in $W(Q,R).$
\end{thm}

Proof. We prove Theorem \eqref{Siroj2} using the method of contracting mappings \cite{L9}. Let be $\Im _{\varepsilon } -$the operator corresponding to the differential expression \eqref{GrindEQ__21_} and conditions \eqref{GrindEQ__22_}-\eqref{GrindEQ__24_}. Let us denote the formal inverse operator by $\Im _{\varepsilon }^{-1} $. The following operator is defined in space $W(Q,R)$:
$$u_{k,\varepsilon }^{(l)} =\Im _{\varepsilon }^{-1} F(u_{k,\varepsilon }^{(l-1)} )\equiv Pu_{k,\varepsilon }^{(l-1)} .$$
1. Let us show that the operator $P$  maps space $W(Q,R) $ into itself.

Let $\left\{\, u_{k,\varepsilon }^{(l-1)} \, \right\}\in W(Q,R)$, then to a solution of  problem \eqref{GrindEQ__21_}-\eqref{GrindEQ__24_}, the assertion of Lemma \eqref{Sir1} is true, i.e. estimate II is valid ). It follows that for any $l=1,2,3...$we obtain $\left\{\, u_{k,\varepsilon }^{(l)} \, \right\}\in W(Q,R)$.

Thus,$P:W(Q\, ,R)\to W(Q,R).$

2. Let us show that $P$  is a contraction operator.

Let $\left\{\, u_{k,\varepsilon }^{(l)} \, \right\},\left\{\, u_{k,\varepsilon }^{(l-1)} \, \right\}\in W(Q,R).$ Consider  new a function $\left\{\, \vartheta _{k,\varepsilon }^{(l)} \, \right\}=\left\{\, u_{k,\varepsilon }^{(l)} \, \right\}-\left\{\, u_{k,\varepsilon }^{(l-1)} \, \right\}$; the assertion of Lemma \eqref{Sir2} is true for it, i.e. the estimate IV is valid , i.e.
$$IV) \quad \; \frac{\varepsilon }{\delta } \left\langle \frac{\partial \Delta ^{2} \vartheta _{k,\varepsilon }^{(l)} }{\partial t} \right\rangle _{0}^{2} +\left\langle \vartheta _{k,\varepsilon }^{(l)} \right\rangle _{4}^{2} \le A_{2} q^{(l)} .$$

Thus, $P$ is a contraction operator. Now, by well-known principle of contraction mappings, problem \eqref{GrindEQ__21_}-\eqref{GrindEQ__24_} has a unique solution belonging to the space $W(Q,R),$ for $\varepsilon >0$. Hence we have $u_{k,\varepsilon }^{(l)} \to u_{k,\varepsilon } $ at $l\to \infty $ \cite{L6}-\cite{L8},\cite{L11},\cite{L12}.

\section{A family of loaded mixed type fourth order differential equations of the second kind}

Now we will prove the unique solvability of problem \eqref{GrindEQ__8_}-\eqref{GrindEQ__11_}. The family of loaded fifth order differential equations \eqref{GrindEQ__21_} with conditions \eqref{GrindEQ__22_}-\eqref{GrindEQ__24_} is used as a `` $\varepsilon -$regularizing'' equation for equation \eqref{GrindEQ__8_} with conditions \eqref{GrindEQ__9_}-\eqref{GrindEQ__11_}  \cite{L3},\cite{L4}-\cite{L12}.

Let, $\left\{\, u_{k,\varepsilon } \, \right\}\in W(Q,R)$ for a fixed value, $\varepsilon >0$ be a unique solution to problem \eqref{GrindEQ__21_}-\eqref{GrindEQ__25_}. Then for $\varepsilon >0$ inequality IV is true. By the weak compactness theorem \cite{L13}-\cite{L24}, from a bounded sequence $\left\{\, u_{k,\varepsilon } \, \right\}$ one can extract a weakly convergent subsequence of the function $\left\{\, u_{k,\varepsilon _{j} } \, \right\} $ such that $u_{k,\varepsilon _{j} } \to u_{k} $ is weak $W(Q,R)$ \cite{L16}. Let us show that the limiting function $u_{k} (x,t) $ satisfies equation \eqref{GrindEQ__8_} almost everywhere in $W(Q,{\mathbb R})$.
 \\ Indeed, since the subsequence $\left\{\, u_{k,\varepsilon _{j} } \, \right\}$ weakly converges in $W(Q,{\mathbb R})$, and operator $\Im -$is linear, we obtain
\begin{equation} \label{GrindEQ__45_}
\begin{array}{l}{\Im u_{k}-F(u_{k})=\Im u_{k}-F(u_{k,\varepsilon _{j} })-[F(u_{k})-F(u_{k,\varepsilon _{j}})]=\varepsilon_{j}\frac{\partial \Delta ^{2} u_{k,\varepsilon_{j}}}{\partial t} +}\\
{+L_{0}(u_{k}-u_{k,\varepsilon_{j}})+\mu_{k}^{4}(u_{k}-u_{k,\varepsilon_{j}}  )-[F(u_{k} )-F(u_{k,\varepsilon_{j}})].} \end{array}
\end{equation}

Passing to the limit in \eqref{GrindEQ__45_} at $\varepsilon _{j} \to 0$, we obtain. $\Im u_{k} =F(u_{k} ).$ Hence the function $u_{k} (x,t)$ will be the unique solution to problem \eqref{GrindEQ__8_}-\eqref{GrindEQ__11_} from $W(Q,R)$ \cite{L5}-\cite{L11}.

\textbf{\textit{This proves Theorem \eqref{Siroj2}. Now let's prove Theorem \eqref{Siroj1}.}}

Since all the conditions of Theorem \eqref{Siroj1} are satisfied, using the Parseval -- Steklov equalities \cite{L11}-\cite{L13} to a solution of  problem \eqref{GrindEQ__8_}-\eqref{GrindEQ__11_}, we obtain a solution to problem \eqref{GrindEQ__1_}-\eqref{GrindEQ__7_} from the specified class $U.$

\begin{remark}The inverse problem for mixed-type equations of the second kind of high even order is studied similarly.\end{remark}

\hfill$\Box$ \
\section*{Acknowledgements}
 The authors are grateful to  Artyushin A.N. for discussions of these results. The authors acknowledge financial support from the  Ministry of Innovative Development of the Republic of Uzbekistan, Grant No F-FA-2021-424.



\begin{thebibliography}{99}
 \normalsize

\bibitem{L1}
Anikanov~Yu.E. \textit{Some methods for studying multidimensional inverse problems for differential equations.} Novosibirsk: Nauka, (1978).

\bibitem{L2}
Berezinsky~Yu.M., \textit{Expansion in eigenfunctions of self-adjoint operators.} Kyiv: Nauk.dumka, (1965).

\bibitem{L3}
Bitsadze A.V.  \textit{Ill-posedness of the Dirichlet problem for equations of mixed type.} DAN SSSR. T.~122, No.~4. ---167--170. (1953).

\bibitem{L4}
Bubnov~B.A. \textit{On the  solvability of multidimensional inverse problems for parabolic and hyperbolic equations.} Novosibirsk . Preprint No.~713, Computing Center SO USSR Academy of Sciences. (1987).

\vskip4pt\bibitem{L5}
Vragov~V.N. \textit{Boundary value problems for non-classical equations of mathematical physics.} Novosibirsk: NSU, (1983).

\bibitem{L6}
Vragov~V.N. \textit{
On the formulation and solvability of boundary value problems for equations of mixed-composite type.} Mathematical analysis and related issues of mathematics. Novosibirsk: IM SO AN USSR,  ---5--13. (1978).

\bibitem{L7}
Egorov~I.E. \textit{Nonclassical equations of high order mathematical physics.} Novosibirsk, (1995).

\bibitem{L8}
Dezin~A.A. \textit{General questions of the theory of boundary value problems.} M: Nauka, (1980).

\bibitem{L9}
Dzhamalov~S.Z. \textit{On the correctness of a nonlocal boundary value problem with constant coefficients for a mixed type equation of the second kind of second order in space.} Mat. NEFU notes, No.~4, ---17--28. (2017).

\bibitem{L10}
Dzhamalov S.Z. \textit{On the smoothness of a nonlocal boundary value problem for a multidimensional equation of mixed type of the second kind in space.} Journal of the Middle Volga Mat Society. Vol.~21, No.~1, ---24--33.(2019).

\bibitem{L11}
Dzhamalov~S.Z. \textit{Nonlocal boundary value and inverse problems for equations of mixed type.} Monograph. Tashkent. (2021).

\bibitem{L12}
Dzhamalov~S.Z., Ashurov~R.R. \textit{On a linear inverse problem for a multidimensional equation of mixed type of the second kind, second order} Differential equations. T.55. No.~1, ---34--44. (2019).

\bibitem{L13}
Dzhamalov~S.Z., Ashurov~R.R. \textit{On a linear inverse problem for a multidimensional equation of mixed type of the first kind, second order.} Russ. Mathematics. No.~6, ---1--12. (2019).

\bibitem{L14}
Dzhamalov~S.Z., Pyatkov~S.G. \textit{Some classes of inverse problems for second order mixed type equations.}
Mat. NEFU notes, T.~25, No.~4, ---3--15.(2018).

\bibitem{L15}
Dzhamalov~S.Z., Pyatkov~S.G. \textit{On some classes of boundary value problems for multidimensional equations of mixed type of high order.}
Siberian Mathematical Journal, T.~61, No.~4, ---777--795.(2020).

\bibitem{L16}
Ladyzhenskaya~O.A. \textit{Boundary value problems of mathematical physics.} M. (1973).

\bibitem{L17}
Lavrentyev~M.M., Romanov~V.G., Vasiliev~V.G. \textit{Multidimensional inverse problems for differential equations.} Novosibirsk: Nauka, (1969).

\bibitem{L18}
Lyons~J.L., Magenes~E.  \textit{Inhomogeneous boundary value problems and their applications.}M.Mir.(1971).

\bibitem{L19}
Hermander Lars \textit{Linear differential operators with partial derivatives.} Publisher: "Mir". Moscow. (1965).

\bibitem{L20}
Nakhushev A.M. \textit{Loaded Equations and their Application. }Nauka, Moscow, (2012) [in Russian].

\bibitem{L21}
Megrabov~A.G. \textit{Some inverse problems for mixed type equations.} DAN SSSR, vol.~234, ---305--307. (1977).

\bibitem{L22}
Sobolev~S.L. \textit{Some applications of functional analysis in mathematical physics.} Moscow: Nauka, (1988).

\bibitem{L23}
Sabitov~K.B., Martemyanova~N.V. \textit{Nonlocal inverse problem for a mixed type equation. }
Russ. Mathematics. No.~2, ---71--85. (2011).

\bibitem{L24}
Kabanikhin S.I. \textit{Inverse and ill-posed problems. }Novosibirsk, Siberian Scientific Publishing House, (2009).

\bibitem{L25}
Kozhanov A.I. \textit{Nonlinear loaded equations and inverse problems.}Computational Mathematics and Mathematical Physics.  T.44, No.4. ---694--716. (2004).

\bibitem{L26}
V.A.Trenogin \textit{ Functional Analysis [Russian].} Moscow: Nauka, (1980).

\bibitem{L27}
Chueshev A.V. \textit{On one linear equation of mixed type of high order.} Siberian Mathematical Journal, vol.43. No.2, ---454--472.(2002).

\bibitem{L28}
 S.Z. Dzhamalov, R.R. Ashurov, A.I. Kozhanov. \textit{Linear inverse problem for 3-dimensional chaplygin equation with semi-nonlocal boundary conditions in a prismatic unbounded domain.} Journal of Mathematical Sciences, Vol. 274, No. 2, August, ---186--200.(2023).

\bibitem{L29}
Dzhamalov S., Khalkhadjayev B., Yusupov Sh. \textit{About the unique solvability of a semi-non-local boundary-value problem for a mixed type equation of the second kind of the fourth order.}
 Bull. Inst. Math., Vol.7. No.2, ---35--41. (2024).

\end{thebibliography}
\end{document}